\documentclass[12pt,twoside]{article}
\usepackage[all]{xy}

\date{}
\begin{document}

\centerline{}

\centerline{}

\centerline{}

\centerline {{\bf  ON k(D)-BLOCKS}}

\centerline{}

\centerline{\Large{\bf}}

\centerline{}

\centerline{\bf {Ahmad M. Alghamdi}}

\centerline{Department of Mathematical Sciences}

\centerline{Faculty of applied Sciences}

\centerline{Umm Alqura University}
\centerline{P.O. Box 14035, Makkah 21955, Saudi Arabia.}
\centerline{e-mail: amghamdi@uqu.edu.sa}

\centerline{}

\newtheorem{Theorem}{\quad Theorem}[section]

\newtheorem{Definition}[Theorem]{\quad Definition}

\newtheorem{Corollary}[Theorem]{\quad Corollary}
\newtheorem{Proposition}[Theorem]{\quad Proposition}

\newtheorem{Lemma}[Theorem]{\quad Lemma}
\newtheorem{Remark}[Theorem]{\quad Remark}
\newtheorem{Example}[Theorem]{\quad Example}
\newtheorem{Exa}[Theorem]{\quad Example}

\begin{abstract}
The objective of   this research paper  is to study the relationship between a block of a finite group and a defect group of such block. We  define a new notion which is called a strongly $k(D)$-block and give  a necessary and sufficient condition of  a block with a cyclic defect group to be  a $k(D)$-block in term of its inertial index.
\end{abstract}



\section{\textbf{Introduction}}
Let $p$ be a prime number and $B$  a $p$-block of a finite group $G$ with a defect group $D$ of order $p^{d}$.
In \cite{1}, we proved   that   $k(B)\equiv  k(D) \ Mod (p)$, where $k(B)$ is the number of the ordinary irreducible characters belonging to  $B$ and $k(D)$ is the number of the ordinary irreducible characters of  the defect group  $D$, which is an extra special $p$-group of order $p^{3}$ and exponent $p$, for an odd prime number $p$. In other words, we  relate  the number of the ordinary irreducible characters of a finite group which belong to a certain block and the number of the conjugacy classes of the defect group of that block under consideration.

This result led  us to  think about numerical relationships between a $p$-block and its defect group.
In the present  work, we are free from any condition about the prime number $p$. A question about the existence of a function on the natural numbers which relates some block-invariants under consideration is known as Brauer's $21$-problem (see \cite{22, 3}). In fact, Brauer  asked whether it is the case that $k(B) \leq |D|$ in general. This question is known as Brauer $k(B)$-conjecture.

 For a bound of $k(B)$, it is well known that $k(B)\leq p^{2d-2}$ for $d\geq 3$ and $k(B)\leq p^{d}$ for $d\in \{0,1,2\}$. See \cite{22}, \cite{3, 4}, \cite{knoor} and \cite{kulrob, kul} for more details and discussions in this direction.

 However, we have arisen a question about which blocks and which conditions  ensure  the equality $k(B)=k(D)$ as well as ensure the congruency $k(B)\equiv  k(D) \ Mod (p)$. We have studied some general cases as well as some examples for small $G$. Then we try to characterize such blocks which have  cyclic defect groups in the term of the order of the inertial subgroups.

However, as far as we know, we have not seen  a similar relation in the literature.  In fact, most of the examples which we have already considered satisfy the equality $k(B)=k(D)$. However, for $G=A_{5}$ with $p=5$, we find that $k(B_{0}(G))=4$, where $B_{0}(G)$ is the principal $5$-block of $G$. But $k(D)=5$, where $D\in Syl_{5}(G)$. Since $A_{5}$ is a simple group and some  unusual properties arise from such group, we think that there are some classes in group theory in such a way, $p$-blocks  satisfy either the equality or the congruency relation. We call those groups which do not satisfy our congruency relation exotic groups.

 \begin{Definition}\label{100}
 Let $p$ be a prime number, $G$ a finite group and $B$ a $p$-block of $G$ with defect group $D$. Write $k(B)$ to mean the number of ordinary irreducible characters of $B$ and write $k(D)$ to mean the number of ordinary irreducible characters of $D$. We call $B$  a strongly  $k(D)$-block if $k(B)=k(D)$.
 \end{Definition}
 We can consider an equality mod the prime number $p$ as in the following definition.
 \begin{Definition}\label{101}
 Let $p$ be a prime number, $G$ a finite group and $B$ a $p$-block of $G$ with defect group $D$. Write $k(B)$ to mean the number of ordinary irreducible characters of $B$ and write $k(D)$ to mean the number of ordinary irreducible characters of $D$. We call $B$  a   $k(D)$-block if $k(B)\equiv k(D) Mod (p)$.
 \end{Definition}
 It is clear that a strongly $k(D)$-block is a $k(D)$-block. However, we shall see in Example \ref{sem} some   $k(D)$-blocks which are not strongly $k(D)$-blocks.

 Our main concern is to study finite groups and their blocks which satisfy Definitions \ref{100} and \ref{101}. Note that it is well known that $k(D)$ is the number of the conjugacy classes of $D$. It is well known that blocks with cyclic defect groups are well understood. This theory is  rich and has many applications. So, we shall start by doing some sort of characterization of strongly $k(D)$-blocks with cyclic defect groups. Our main tool is Dade theorem for the number of irreducible characters of a block with a cyclic defect group; ( see \cite{dade} and \cite[Page 420]{lary}).

 In the end of the paper, we use the computations and the  results in \cite{kulsamb, kessar} to see that such phenomena do occur quite often in block theory.

\subsection{Examples of strongly  $k(D)$-blocks, $k(D)$-blocks  and non $k(D)$-blocks:}
\label{exam}
We shall start with some examples which illustrate the phenomenon  of $k(D)$-blocks.
\begin{Exa}\label{sem}
For $G=S_{n}$; the symmetric group of $n$ letters and $n\in\{2,3, 4, 5, 6, 7, 8\}$, it happens  that $k(B_{0}(G))=k(D)$, for $D\in Syl_{p}(G)$ and $p\in \{2, 3, 5, 7\}$. However, for any  prime number  $p\geq 5$, the defect group of the principal $p$-block of the symmetric group $S_{2p}$ is an abelain $p$-group of order $p^{2}$ and $k(B_{0})=\frac{p^{2}+3p}{2}$. In this situation, we obtain  $k(D)$-blocks which are not strongly $k(D)$-blocks. A similar conclusion  holds for the principal $p$-block of the symmetric group $S_{3p}$, with $p\geq 7$. However, when $p=3$ and $G=S_{9}$, then $22=k(B_{0}(S_{9}))\neq k(D)=17$, where $D\in Syl_{3}(S_{9})$.
\end{Exa}
\begin{Exa}
Let $G$ be the dihedral group of order $8$. It has a unique $2$-block with five ordinary irreducible characters which obviously coincide with the number of the conjugacy classes of $G$.
\end{Exa}

\begin{Exa}
Let $G$ be the alternating group $A_{4}$. Then for $p=2$, $G$ has a unique $2$-block with four ordinary irreducible characters which is the same as the number of the conjugacy classes of
 the defect group. For $p=3$, we see that $G$ has one $3$-block of defect zero and the principal $3$-block with three ordinary irreducible characters with the same number as the number of the conjugacy classes of a Sylow $3$-subgroup of $G$.
\end{Exa}
\begin{Exa}
For $G=SL(2, 3)$; the special linear group,  we have for $p=2$, the principal $2$-block has the quaternion group $Q_{8}$ as a defect group and indeed, $k(B_{0}(SL(2, 3)))=7$ and $5=k(Q_{8})$. Then  $k(B_{0}(SL(2, 3)))$ is a $k(D)$-block.
\end{Exa}

\begin{Exa}
Now, we have faced the first example which does not obey  our speculation. It is the first non abelian simple group; $A_{5}$. Although, for $p\in \{2,3\}$, $k(B_{0}(A_{5}))=k(D)$, where $D\in Syl_{p}(A_{5})$, we observed  that for $p=5$, $k(B_{0}(A_{5}))=4\neq k(D)=5$, where, $D\in Syl_{5}(A_{5})$. Hence, $A_{5}$ is the first exotic group in this sense.   The same obstacle we have faced for the  group $G=GL(3, 2)$, since $k(B_{0}(G))\neq k(D),$ when $D\in Syl_{7}(G)$
\end{Exa}
\begin{Exa}
The principal $3$-block for $G=PSL(3, 3)$; the projective special linear group,  satisfies $11=k(B_{0}(G))=k(D)=11,$ where $D$ is an extra special $3$-group of order $27$ and exponent $3$.
\end{Exa}

 \subsection{General cases  for the notion of $k(D)$-blocks}
 \label{gcase}
 \begin{enumerate}
\item Let  $p$ be a prime number and $G$  a finite group. Assume that the prime number $p$ does not divide the order of the group $G$. Then each block of $G$ has  defect zero, (see Theorem 6.29 Page  247 in  \cite{2} ). Hence such block  is a strongly  $k(D)$-block.
\item It is well known,  see Problem 13 in Chapter 5 page 389 in \cite{2} that if $G=DO_{p^{'}}(G)$ is a $p$-nilpotent group  with a  Sylow $p$-subgroup $D$ and the maximal normal  $p^{'}$-subgroup of $G$ then  $k(B)=k(D)$. Certainly,  a nilpotent $p$- block is  a strongly  $k(D)$-block.
\item We know that if $G$ is a $p$-group then it has a unique $p$-block, namely, the principal block and such block is an strongly  $k(D)$-block.
\item For $p$-blocks with dihedral defect groups, we have $k(B)=k(D)$. Hence such blocks are examples of strongly  $k(D)$-blocks.
\end{enumerate}

 \section{ $k(D)$-blocks with cyclic defect groups}
 \label{gca}
 In this section, we discuss $p$-blocks with cyclic defect groups.
Recall that a  root of  a $p$-block $B$ of a finite group $G$ with defect group $D$ is a $p$-block $b$ of the subgroup $DC_{G}(D)$ such that $b^{G}=B$ (see  Chapter 5 Page 348 in \cite{2}),  where $C_{G}(D)$ is the centralizer subgroup of $D$ in $G$. Now for the root  $b$, we define the inertial index of $B$ to be the natural number $e:=e(B)=[I_{G}(b): DC_{G}(D)]$, where $I_{G}(b)=\{g\in G: b^{g}=b\}$. It is clear that $I_{G}(b)$ is a subgroup of $G$ which contains $DC_{G}(D)$ and the index $[I_{G}(b): DC_{G}(D)]$ is well-defined. The number $e$ above is crucial to investigate some fundamental results in block theory.

Let us restate the following well known result which was established by Dade regarding the number of irreducible characters in a block with a cyclic defect group. For more detail, the reader can see the proof and other constructions in \cite{dade, nav}.
\begin{Lemma}\label{dade}
Let  $B$ be a $p$-block of a  finite group $G$ with a cyclic defect group $D$ of order $p^{d}$. Then $B$ has $e+\frac{p^{d}-1}{e}$ ordinary irreducible characters, where $|D|=p^{d}$.
\end{Lemma}

With the above notation, we characterize strongly $k(D)$-blocks in the term of the inertial index for  blocks with cyclic defect groups. Also, we believe that it is worth looking for some positive theorems regarding the notion of $k(D)$-blocks.
\begin{Theorem}\label{001}
Let  $B$ be a $p$-block of a finite group $G$ with a cyclic defect group $D$ of order $p^{d}$. Then $B$  is a strongly $k(D)$-block if and only if $e=1$ or ($e=p-1$ and $d=1$).
\end{Theorem}
Proof: Assuming that $B$ is a strongly $k(D)$-block and using  Lemma \ref{dade}, we can write $p^{d}=e+\frac{p^{d}-1}{e}.$ Then we have $e^{2}-p^{d}e+(p^{d}-1)=0.$ Letting $e$  be the variable, we see that the only solution we have is that $e=1$ or $e=p^{d}-1$. The result follows as $e$ divides $p-1$. The converse is clear and the main result follows.

\begin{Remark}
We  get an analogue  result of Theorem \ref{001} for $k(D)$-blocks with cyclic defect groups, by solving  the congruency equation  $p^{d}\equiv_{p} e+\frac{p^{d}-1}{e}.$
\end{Remark}

\section{The interplay with fundamental  results}
There are fundamental progress in solving Brauer problems. We recast the following result which is due to R. Kessar and G. Malle \cite[HZC1]{kessar}. This result can be used to see an strongly block $B$ with abelian defect group  $D$ of order $p^{d}$  as such that $k_{0}(B)=k(D)$, where $k_{0}(B)$ is the number of ordinary irreducible characters of height zero belonging to $B$.
\begin{Lemma}
Let $G$ be a finite group, $B$ be a $p$-block of $G$ with defect group $D$. If $D$ is abelian then every ordinary irreducible character of $B$ has height zero.
\end{Lemma}

Let us conclude this paper by mentioning the following lemma in such a way that we rely on the computation  in \cite[Proposition 2.1]{kulsamb} by B. Kulshammer and D. Sambale.  These computations guarantee that the phenomena of strongly $k(D)$-block occur quite often in the theory of blocks.

\begin{Lemma}\label{rep}
Let $G$ be a finite group, $B$ be a $p$-block of $G$ with a defect group $D$. If $D$ is an elementary abelian of order 16 then  $B$ is a $k(D)$-block.
\end{Lemma}
In fact, Lemma  \ref{rep} can be replaced by the following much stronger result.

\begin{Theorem}\label{0123}
Let   $G$ be a finite group and $p=2$ or $3$. Then each $p$-block of $G$ with abelian defect group is a  $k(D)$-block.
\end{Theorem}


\bigskip



\begin{thebibliography}{99}
\bibitem{1}  A. Alghamdi, The ordinary weight conjecture and Dade's
projective conjecture for $p$-blocks with an extra special defect
group, Ph.D. Thesis, University of Birmingham, 2004.

\bibitem{22} R. Brauer and W. Feit, \emph{On the number of irreducible characters of finite groups in a given block,} Proc. Nat. Acad. Sci. U. S. A. 45 (1959), 361-365.
    127--131.

    \bibitem{dade} E. Dade, \emph{Counting characters in blocks with cyclic defects groups,} J. Algebra 186 (3), 934-969, (1996).

        \bibitem{dade66} E. Dade, \emph{Blocks With Cyclic Defect Groups,} Annals of Mathematics, Second Series, Vol. 84, No. 1 (Jul., 1966), pp. 20-48.


        \bibitem{lary} L. Dornhoff, \emph{Group Representation Theory,} Part B:  Modular representation theory. Marcel Dekker  Inc., New York,  (1972).

\bibitem{knoor} R. Knorr, \emph{On the number of characters in a $p$-block of a $p$-solvable group,} Illinois J. Math.  28, 181-210 (1984).

\bibitem{kulrob} B. Kulshammer and G. Robinson, \emph{Alperin-Makey Implies Brauer's Problem 21,} J. Algebra 180, 208--210 (1996).

\bibitem{kul} B. Kulshammer, \emph{Modular representations of finite groups: conjectures and examples,}  Jena,  (1996).
\bibitem{kulsamb} B. Kulshammer and B. Sambale, \emph{The 2-blocks of defect 4,}  Journal: Represent. Theory 17 (2013), 226-236, Published electronically: May 2, 2013.

\bibitem{kessar} R. Kessar and G. Malle, \emph{Quasi-isoloated blocks and height zero Conjecture,} J. Algebra 147, 450--455 (1992).


\bibitem{2} H. Nagao and Y. Tsushima, \emph{Representation of finite groups, } Academic Press Inc., Boston, MA,  Translated from Japanese, (1989).

    \bibitem{nav} G. Navarro, \emph{Characters and blocks of finite groups,} Volume 250 of London Mathematical Society Lecture Notes Series. Cambridge University Press, Cambridge 1998.

\bibitem{3} G. R. Robinson, \emph{Some open conjectures on representation theory,} J.
 Ohio State Univ. Math. Res. Inst. Publ.,de Gruyter,
   Berlin,1997,  Volume \textbf{6}, Representation theory of finite groups (Columbus, OH, 1995),
    127--131.

    \bibitem{4} G. R. Robinson, \emph{On Brauer's $k(B)$ Problem,} J. Algebra 147, 450-455 (1992).
    127--131.
\end{thebibliography}
\end{document}